\documentclass[12pt]{amsart}
\usepackage{amsmath, amssymb, euscript, graphicx, hyperref}

\newtheorem{thm}{Theorem}[section]

\newtheorem{prop}[thm]{Proposition}
\theoremstyle{definition}
\newtheorem{defn}[thm]{Definition}

\newcommand{\blackboard}[1]{\ensuremath{\mathbb{#1}}}

\newcommand{\N}{\blackboard{N}}

\newcommand{\Z}{\blackboard{Z}}

\newcommand{\Q}{\blackboard{Q}}
\newcommand{\R}{\blackboard{R}}
\newcommand{\C}{\blackboard{C}}

\begin{document}

\author{Azer Akhmedov}
\address{Azer Akhmedov,
Department of Mathematics,
North Dakota State University,
PO Box 6050,
Fargo, ND, 58108-6050}
\email{azer.akhmedov@ndsu.edu}

\title{Irreducible discrete subgroups in products of simple Lie groups}

\maketitle


\

  \medskip

\maketitle

\begin{center} {\small ABSTRACT: We produce an example of an irreducible discrete subgroup in the product $SL(2,\R)\times SL(2,\R)$ that is not a lattice. This answers a question asked in \cite{FMvL}.}\footnote{In the forthcoming update of this paper, we will extend the result (i.e. Theorem \ref{thm:main}) to all isotypic products of simple Lie groups with at least two factors and without a compact factor. The proof of this more general result uses the main idea of the current version.} \end{center}

\section{Introduction}

 We are motivated by the following question of Fisher-Mj-van Limbeek (see Question 1.6 in \cite{FMvL}): 

\medskip 

  {\em Question 1.} Let $G_1$ and $G_2$ be semisimple groups over local fields and
$\Gamma \leq G_1\times G_2$ be a discrete subgroup with both projections dense. Is $\Gamma $ in fact an irreducible lattice in the product $G_1\times G_2$?

\medskip 

 According to \cite{FMvL}, Yves Benoit first informally posed this question as far back as 2003 or 2004. \cite{FMvL} and \cite{BFMvL} discuss some very interesting motivations for this question relating it also to the following question of Greenberg-Shalom: 

  \medskip

   {\em Question 2.} Let $G$ be a semisimple Lie group with finite center and without compact factors. Suppose $\Gamma \leq G$ is a discrete, Zariski-dense subgroup of $G$ whose commensurator $\Delta \leq G$ is dense. Is $\Gamma $ an arithmetic lattice in $G$?

\medskip 

   In its own turn, Question 2 is strongly motivated by Margulis-Zimmer Conjecture (see \cite{SW}
and Conjecture 1.4 of \cite{FMvL}). Let us note that the object $\Delta $ in Question 2 is an important characterizing object; by a landmark theorem of Margulis, $\Delta $ detects arithmeticity of lattices (in real or $p$-adic semisimple Lie groups with finite center and without compact factors; see \cite{M} or Theorem 1.1. in \cite{FMvL}).    

\medskip 

 We provide a negative answer to Question 1 by constructing a discrete free subgroup in a product with dense projections. Let us recall that an irreducible lattice $\Gamma $ in a higher rank semi-simple Lie group without a compact factor has the following property: 

\medskip 

 (i) $\Gamma $ contains a copy of $\Z ^2$ (see \cite{PR}). 

\medskip 

 By property (i), $\Gamma $ cannot be free. Non-freeness is a weak (still a meaningful) property of higher rank irreducible lattices, but it is the easiest (that we found) to use to produce an example that needed for Question 1.  

\medskip 

  Thus our aim will be to construct free discrete subgroups (with dense projections).  

  \medskip 
  
  Both freeness and discreteness of subgroups can be difficult to establish in various contexts/environments. For connected Lie groups, there are elementary open questions in this area even for $SL(2,\R)$ \cite{B}, \cite{KK}, \cite{LU}. For the group $\mathrm{Diff}_{+}(I)$, the $C_0$-discreteness has been studied in \cite{A1} where a characterization of such groups have been presented in $C^{1+\epsilon }$ regularity. A more complete characterization of such groups has been presented in \cite{A2} and \cite{A3}.   In \cite{A4}, the $C^1$-discreteness question has been discussed and some elementary open questions have been raised. We also refer the reader to \cite{A5} which directly studies free and discrete subgroups of $\mathrm{Diff}_{+}(I)$. It is interesting that the $\mathrm{Homeo}_{+}(I)$ and $\mathrm{Diff}_{+}(I)$ environments provide other tools (not discussed in this paper) for establishing freeness and discreteness of certain subgroups. We refer the reader to a series of remarkable papers \cite{AGM}, \cite{C}, \cite{C1},  \cite{CG}, \cite{CG2}, \cite{W} which establish existence of free subgroups. 

  \medskip 

  Our main result is the following theorem.

  \begin{thm} \label{thm:main} There exists a free discrete subgroup $\Gamma \leq SL(2,\R )\times SL(2,\R )$ such that the projection of $\Gamma $ to each factor is dense. 
  \end{thm} 

  \medskip 

  In our proof, we first consider the case of $ SL(2,\R )\times SL(2,\C )$ and prove the following theorem.

   \begin{thm} \label{thm:main2} There exists a free discrete subgroup $\Gamma \leq SL(2,\R )\times SL(2,\C )$ such that the projection of $\Gamma $ to each factor is dense. 
  \end{thm} 

   \medskip 
   
  As pointed out in \cite{FMvL}, it is easy to produce an example with a dense projection in one of the factors.  Let us also recall a classical fact that the group $SL(2,\Z [\sqrt{2}]$ is an irreducible lattice in $SL(2,\R)\times S(2,\R )$ by the faithful representation $\rho : SL(2,\Z [\sqrt{2}])\to SL(2,\R)\times S(2,\R )$ given by $\rho (A) = (A, \sigma (A)), A\in SL(2,\Z [\sqrt{2}])$ where $\sigma :SL(2,\Z [\sqrt{2}])\to SL(2,\Z [\sqrt{2}])$ is the Galois isomorphism obtained from the Galois automorphism $\sigma : \Z [\sqrt{2}]\to \Z [\sqrt{2}]$ of the ring $\Z [\sqrt{2}]$ defined as $\sigma (m+n\sqrt{2}) = m-n\sqrt{2}, m. n\in \Z$. The discreteness of $\rho (SL(2,\Z [\sqrt{2}]))$ comes from the fact that in it, if $\rho (A)$ converges to identity in one factor, then it escapes to infinity in the other. This phenomenon also causes difficulty (among other issues) in attempts to construct a straightforward example for the claim of Theorem \ref{thm:main}. Let us also recall that by the main result of \cite{BG}, every dense group in a semi-simple Lie group contains a dense free subgroup; this result seems somewhat relevant here, but in trying to apply it (or the idea of it), one has to fight this time to preserve the discreteness of a subgroup in the product. Question 1 is indeed at a very interesting conjuncture of tensions among freeness, discreteness, and denseness. Another manifestation of this lies in the fact that, to establish freeness, it is more suitable to use hyperbolic or parabolic elements, whereas for denseness, elliptic elements are more efficient. 

  \medskip 

   Given a subgroup $\Gamma $ in a product $G_1\times G_2 \times \ldots \times G_n$ of simple non-compact Lie groups with $n\geq 2$, we call $\Gamma $ {\em irreducible} if for all $1\leq i\leq n$, $\pi _i(\Gamma )$ is dense in $G_i$ where $\pi _i:G\to G_i$ is the projection onto the $i$-th factor. Thus, Theorem \ref{thm:main} establishes the existence of a discrete irreducible subgroup in the product  $SL(2,\R )\times SL(2,\R )$ which is not a lattice.  In the forthcoming update of this paper, we will extend the result to all products $G = G_1\times G_2 \times \ldots \times G_n$ of simple Lie groups with $n\geq 2$ where the group $G$ is {\em isotypic} (i.e. all simple factors of $G_{\mathbb{C}}$ are isogenous to each other) and has no compact factors. Recall that by a result of Margulis, if $G$ has no compact factors and admits an irreducble lattice, then it is isotypic.  The converse (even without the assumption about compact factors) also holds, cf. \cite{Mo}. 

\bigskip 

  {\em Acknowledgment:} I am very thankful to David Fisher, Tsachik Gelander, and Yehuda Shalom for being interested in this work and for pointing out a serious error in the earlier version of this paper. I am also grateful to Emmanuel Breuillard for very useful discussions on the previous version of this paper. The author has used ChatGPT to perform some auxiliary computations (not directly used in the paper). 

\vspace{1cm}

     \section{Matrices with dominant eigenvalues} 

     \medskip 
     
     In this section, we will briefly review the tools used in the proof of Tits Alternative for linear groups. Our terminology is somewhat different from the one used in \cite{D}. In particular, we will restrict ourselves to the real case, but the discussions can be generalized to any locally compact normed field (in particular, to any local field) as it is done in \cite{D}. 

\medskip 

      \begin{defn} An eigenvalue $\lambda $ of a matrix $C\in GL(n,\R)$ is called {\em dominant} if $\lambda $ is real, has multiplicity 1 and $|\lambda | > \max \{|\mu |, 1\}$ for any other eigenvalue $\mu $ of $C$. A matrix $C\in GL(n,\R)$ is called {\em hyperbolic-like} if both $C$ and $C^{-1}$ have dominant eigenvalues. If $n=2$, then a hyperbolic-like matrix is just called hyperbolic. 
      \end{defn}

\medskip 

   The group $GL(n, \mathbf{k})$, for any field $\mathbf{k}$, has a standard action on $\mathbb{P}_{\mathbf{k}}^{n-1}$. If $C\in GL(n, \R)$ is hyperbolic-like, then it has unique and distinct {\em attracting and repelling points} $a, b\in \mathbb{P}_{\R}^{n-1}$, and {\em characteristic crosses} $\Pi _a, \Pi _b$ such that for any compact $K_1\subseteq \mathbb{P}_{\R}^{n-1}\backslash \Pi _a, K_2\subseteq \mathbb{P}_{\R}^{n-1}\backslash \Pi _b$ there exist open neighborhoods $U_1, U_2$ of $a, b$ respectively, and a natural $N\geq 1$ such that for all $n\geq N$, $A^n(K_i)\subseteq U_i, 1\leq i\leq 2$. The points $a, b$ are indeed (the class of) eigenvectors of $C$, corresponding to the biggest and smallest eigenvalues; $\Pi _a = \mathbb{P}(V_a), \Pi _b = \mathbb{P}(V_b)$ are projectivizations of subspaces $V_a, V_b\subset \mathbb{R}^n$ of dimensions $n-1$ (so  $\Pi _a, \Pi _b$ are subvarieties of  $\mathbb{P}_{\R}^{n-1}$ of dimension $n-2$; moreover, $a\notin \Pi _a, b\notin \Pi _b$ and $a\in \Pi _b, b\in \Pi _a$). If we let $\Lambda _C$ be the set of eigenvalues, with $\lambda , \mu $ being the biggest and the smallest eigenvalues, then, in the subspaces $V_a, V_b$ are associated with the set of eigenvalues $\Lambda \backslash \{\lambda \}$ and $\Lambda \backslash \{\mu \}$ respectively. In the case when $\Lambda \subset \R$, we have $V_a = Span (\Lambda \backslash \{\lambda \})$ and $V_b = Span (\Lambda \backslash \{\mu \})$.
   
   \medskip

   We will use the notation $A_C, R_C$ for attractive and repelling points of $C$ respectively and write $\mathcal{F}_C = \{A_C, R_C\}$, i.e. $A_C := a, R_C := b$. We also will write $\Pi _C^+ := \Pi _a, \Pi _C^- := \Pi _b$ and $\Pi _C = \Pi _C^+\cup \Pi _C^-$. Let us note that if $C$ is hyperbolic-like, then for all $n\in \Z \backslash \{0\}$, $\mathcal{F}_{C^n} = \mathcal{F}_C$; moreover, if $n > 0$, then $A_{C^n} = A_C, R_{C^n} = R_C$ and if $n<0$, then  $A_{C^n} = R_C, R_{C^n} = A_C$. By a standard ping-pong argument, we will obtain the following proposition.

    \medskip 

    \begin{prop} \label{prop:dominant} Let $n\geq 1$, and $A, B\in GL(n,\R)$ be hyperbolic-like matrices such that $\mathcal{F}_A\cap (\mathcal{F}_B \cup \Pi _B) = \mathcal{F}_B\cap (\mathcal{F}_A \cup \Pi _A) = \emptyset $. Then there exists $N\geq 1$ such that for all $m, k\geq N$, the matrices $A^m, B^k$ generate a discrete free group of rank two.  
    \end{prop}

\medskip 

    The discreteness of the subgroup $\langle A^m, B^k\rangle $ in the above proposition is meant in the natural topology of $GL(n, \R)$. Let us note that, for any hyperbolic-like matrix $C\in GL(n, \R)$, we also have $\mathcal{F}_C\subset \Pi _C$ therefore the statement of Proposition \ref{prop:dominant} can be simplified by observing that $\mathcal{F}_A \cup \Pi _A = \Pi _A$ and $\mathcal{F}_B \cup \Pi _B = \Pi _B$. 

\medskip 

    We also would like to note (recall) the following easier fact which will be used in the sequel as well.

    \begin{prop} \label{prop:noncommuting} Let $n\geq 1$, and $A, B\in GL(n,\R)$ be hyperbolic-like matrices such that $\mathcal{F}_A\cap \mathcal{F}_B = \emptyset $ and $A(\mathcal{F}_B)\cap \mathcal{F}_B = \emptyset $. Then $AB\neq BA$.   
    \end{prop}

    \medskip 

    The results quoted in this section will be applied in a specific case. We would like to state a proposition preparing a setting in which we will deduce the existence of a free subgroup. 

    \begin{prop} \label{thm:lastfactor} Let $\rho :G\to GL(n,\C)$ be a representation of a group $G$ which is a direct sum $\rho _1\oplus \rho _2$ of two representations  $\rho _i:G\to GL(n_i, \C), 1\leq i\leq 2$ with $n_2 = 2$. Let also $a, b\in G$ such that the matrices $\rho _2(a)$ and $\rho _2 (b)$ are hyperbolic without a common eigenvector. Then there exists $N\geq 1$ such that for all $m, n\geq N$, the group $\langle a^m, b^m\rangle $ is a non-Abelian free group.

        \end{prop}

    \vspace{1cm} 
    
 \section{Proof of Theorem \ref{thm:main2} }

  We will build a homomorphism $\mathbf{\Phi } :SL(2,\Z [\sqrt{r}])\to GL(4, \R)$ as a key tool in the proof of Theorem \ref{thm:main}. A crucial property of this homomorphism will lie in the fact that even for elliptic matrices $A\in SL(2,\Z [\sqrt{r}])$, the associated matrix $\mathbf{\Phi }(A)$ in $GL(4,\R)$ can still be hyperbolic-like.    

   \medskip 

   Let $r\geq 2$ be a square-free integer. Any matrix $A\in SL(2,\Z [\sqrt{r}])$ acts on $\Z [\sqrt{r}]^2$. Any $\theta \in \Z [\sqrt{r}]^2$ can be written as  $X+\sqrt{r}Y\in \Z [\sqrt{r}]^2$ with $X =   \begin{bmatrix} 
	x  \\
	z    
	\end{bmatrix} $ and $Y =   \begin{bmatrix} 
	y  \\
	t  \\   
	\end{bmatrix} $. Letting 
 $A =   \begin{bmatrix} 
	a+m\sqrt{r}  & b+n\sqrt{r}   \\
	c+k\sqrt{r} & d+l\sqrt{r} \\   
	\end{bmatrix} $ we obtain that $A\theta =   \begin{bmatrix} 
	ax+rmy+bz+rnt  \\
	cx+rky+dz+rlt    
	\end{bmatrix} +   \sqrt{r} \begin{bmatrix} 
	mx+ay+nz+bt  \\
	kx+cy+lz+dt  \\   
	\end{bmatrix} $.

 \medskip 

  Then we define $\mathbf{\Phi}(A) =    \begin{bmatrix} 
	a & rm & b & rn \\
	m & a & n & b \\
        c & rk & d & rl \\
        k & c & l & d
	\end{bmatrix} $. One can check directly that $\mathbf{\Phi}$ is a monomorphism.

   \medskip 

    Now, we will choose $r = 2$ (we could still work with any square-free $r$). The homomorphism $\Phi $ provides a faithful  representation of the lattice $SL(2,\Z [\sqrt{2}])$. By Margulis Superrigidity Theorem, this representation lifts to a representation of $SL(2,\R )\times SL(2, \R)$. This lift, as pointed out to me by D.Fisher, T.Gelander and Y.Shalom, decomposes into the sum of two standard representations of $SL(2, \R)$. As such, we will not be able to use $\Phi $, to make arrangements satisfying conditions of Proposition \ref{prop:dominant}. \footnote{ Let us also point out that the machine estimates used in the last section of the previous version were not accurate; this lead to an incorrect conclusion there.} This motivates us to consider more sophisticated versions of $\Phi $. 

\medskip 

 First, we will consider the rings $\Z [\sqrt[3]{2}]$ and $\Q [\sqrt[3]{2}]$. Any matrix $A\in SL(2,\Q [\sqrt[3]{2}]$ acts on $\Q [\sqrt[3]{2}]^2$ and any $\theta \in \Q [\sqrt[3]{2}]^2$ can be written as  $X+\sqrt[3]{2}Y+\sqrt[3]{4}Z \in \Q [\sqrt[3]{2}]^2$ with $X =   \begin{bmatrix} 
	x  \\
	u    
	\end{bmatrix} $, $Y =   \begin{bmatrix} 
	y  \\
	v  \\   
	\end{bmatrix} $ and $Z =   \begin{bmatrix} 
	z  \\
	w    
	\end{bmatrix} $ in $\Q^2$.
 
  \medskip 
  
  Letting 
 $A =   \begin{bmatrix} 
	a+m\sqrt[3]{2}+p\sqrt[3]{4}  & b+n\sqrt[3]{2} + q\sqrt[3]{4}   \\
	c+k\sqrt[3]{2} + r\sqrt[3]{4} & d+l\sqrt[3]{2} + s\sqrt[3]{4}\\   
	\end{bmatrix} $ we obtain that $$A\theta =   \begin{bmatrix} 
	ax+bu+2py+2qv+2mz+2nw  \\
	cx+du+2ry+2sv+2kz+2lwlt    
	\end{bmatrix} +   \sqrt[3]{2} \begin{bmatrix} 
	mx+nu+ay+bv+2pz+2qw  \\
	kx+lu+cy+dv+2rz+2sw  \\   
	\end{bmatrix} + $$ \ \ $$\sqrt[3]{4} \begin{bmatrix} 
	px+qu+my+nv+az+bw  \\
	rx+su+ky+lv+cz+dw  \\   
	\end{bmatrix}.$$

\medskip 

 This motivates us to define $\Phi_3 :SL(2, \Q [\sqrt[3]{2}])\to GL(6,\R)$ by letting $$\Psi (A) =    \begin{bmatrix} 
	a & 2p & 2m & b & 2q & 2n \\
	m & a & 2p & n & b & 2q \\
        p & m & a & q & n & b \\
        c & 2r & 2k & d & 2s & 2l \\
        k & c & 2r & l & d & 2s \\
        r & k & c & s & l & d
	\end{bmatrix} .$$ \ As one can check directly, $\Phi_3 $ also turns out to be a monomorphism. 

   \medskip 

   For any integer $\kappa \geq 2$, considering the ring $\Z [\sqrt[\kappa ]{2}]$, we can similarly define the monomorphism $\Phi _{\kappa }: SL(2,\Z [\sqrt[\kappa ]{2}])\to GL(2\kappa ,\R)$ (for $\kappa =2$, we obtain our original map $\Phi $ as $\Phi _2$.) For $\kappa =4$, we obtain the monomorphism $\Psi := \Phi_4 :SL(2, \Z [\sqrt[4]{2}])\to GL(8,\R)$ by letting 

$$\Psi (A) =    \begin{bmatrix} 
	a & 2e & 2p & 2m & b & 2f & 2q & 2n \\
	m & a & 2e & 2p & n & b & 2f & 2q \\
        p & m & a & 2e & q & n & b & 2f \\
        e & p & m & a & f &  q & n & b \\
        c & 2g & 2r & 2k & d & 2h & 2s & 2l   \\     
       k & c & 2g & 2r & l & d & 2h & 2s \\
        r & k & c & 2g & s & l & d & 2h \\
       g & r & k & c & h & s & l & d
	\end{bmatrix} .$$

    Now, we let $$P =  \begin{bmatrix} 
	5+ \beta ^2 - 3\beta -2\beta ^3 & 1 \\
	-1 & 0 
 \end{bmatrix} \ \mathrm{and}  \ Q = \begin{bmatrix} 
	3 + 2\beta ^2 & 1 \\
	-1 & 0 
 \end{bmatrix};$$ here, and for the rest of the paper, $\beta $ will denote the number $\sqrt[4]{2}$. Then $$\Psi (P) =  \begin{bmatrix} 
	 5 & -4  & 2  & -6 & 1 & 0 & 0 & 0 \\
	-3 & 5  & -4 & 2 & 0 & 1 & 0 & 0 \\
         1 & -3 & 5 & -4 & 0 & 0 & 1 & 0 \\
        -2 & 1 & -3 & 5 & 0 & 0 & 0  & 1 \\
         -1 & 0 & 0 & 0 & 0 & 0 & 0 & 0 \\        
         0 & -1 & 0 & 0 & 0 & 0 & 0 & 0 \\
        0 & 0 & -1 & 0 & 0 & 0 & 0 & 0 \\
        0 & 0 & 0 & -1 & 0 & 0 & 0 & 0 
	\end{bmatrix} $$ and $$\Psi (Q) =  \begin{bmatrix} 
	3 & 0 & 4 & 0 & 1 & 0 & 0 & 0 \\
	0 & 3 & 0 & 4 & 0 & 1 & 0 & 0 \\
        2 & 0 & 3 & 0& 0 & 0 & 1 & 0 \\
        0 & 2 & 0 & 0 & 0 & 0 & 0 & 1 \\
        -1 & 0 & 0 & 0 & 0  & 0 & 0 & 0 \\
        0 & -1 & 0 & 0 & 0 & 0 & 0 & 0 \\
        0 & 0 & -1 & 0 & 0 & 0 & 0 & 0 \\
        0 & 0 & 0 & -1 & 0 & 0 & 0 & 0 \\

	\end{bmatrix} .$$

  \medskip 
  
  Notice that the matrix $P$ is elliptic and the matrix $Q$ is hyperbolic. Let $\sigma _k: \Z[\beta ]\to \C, k\in \{0, 1,2,3\}$ be the Galois embeddings defined by $\beta \to \beta \bf{i}^k$ (the map $\sigma _0$ acts as the identity). Then the matrices $\sigma _1(P), \sigma _2(P), \sigma _3(P)$ are all  hyperbolic, whereas $\sigma _0(P) = P$ is elliptic. On the other hand, the matrices $\sigma _1(Q), \sigma _3(Q)$ are elliptic whereas the matrices $Q, \sigma _2(Q)$ are hyperbolic.

\medskip

The representation $\Psi :SL(2, \Z [\sqrt[4]{2}])\to GL(8,\R)$ can be extended by the same formula to $\Psi :SL(2, \Q [\sqrt[4]{2}])\to GL(8,\R)$ with a broader domain which is dense in $SL(2, \R)$. In fact, $\Psi $ can be lifted to the representation of $G := SL(2,\R)\times  SL(2,\R)\times SL(2,\R)\times SL(2,\R)$; it will be the sum of four representations of $SL(2,\R)$ obtained by restricting the lift to each of the four factors of $G$ separately. Notice that by the embedding $x\to (x, \sigma _2(x), \sigma _1(x))$, we can realize $SL(2, \Z [\sqrt[4]{2}])$ as a lattice of $H := SL(2,\R)\times SL(2,\R)\times SL(2,\C)$. Then, by a result of Borel-Harish-Chandra, $\Psi $ can be lifted to a representation of $H$. On the other hand, for the representations of $sl_2\C $, passing to the representations $so(1,3)$ and identifying the complexification of the latter with the complexification of $su_2\oplus su_2$ , we obtain the lift of $\Psi $ to $G$. In more elementary terms, this means that for any real representation of $SL(2,\C)$, the restriction of it to $SL(2,\R)$ is a direct sum of two representations of $SL(2,\R)$, hence it is a representation of $SL(2,\R)\times SL(2,\R)$. However, interestingly, this lift, i.e. the restriction of the lift of $\Psi $ to $G$ is not suitable for our purposes. This is the reason, in the last section we use a special argument to take care of the case of $SL(2,\R)\times SL(2,\R)$ to prove Theorem \ref{thm:main}. For the moment (in this section and in the next section), we restrict our attention to the proof Theorem \ref{thm:main2} and consider the representation $\Psi :H\to GL(8,\R)$.

\medskip 

The lifted representation $\Psi $ can be conjugated to block-diagonal form with each block consisting of size $2\times 2$; for any $A\in SL(2, \Z[\beta ])$, the matrix $\Psi (A)$ will act in each factor as matrices $A, \sigma _1(A), \sigma _2(A)$ and $\sigma _3(A)$ respectively. Hence, the set of eigenvalues of $\Psi (P)$ will be the union of the set of eigenvalues of $P, \sigma _1(P), \sigma _2(P)$ and $\sigma _3(P)$. Similarly, the set of eigenvalues of $\Psi (Q)$ will be the union of the set of eigenvalues of $P, \sigma _1(Q), \sigma _2(Q), \sigma _3(Q)$. Thus, it can be verified directly that the matrix $\Psi (P)$ is hyperbolic-like whereas $\Psi (Q)$ is not; for the latter, we have a double maximal eigenvalue (originating from the factors of $Q$ and $\sigma _2(Q)$).

  \medskip
  
  The eigenvectors of $\mathbf{\Psi }(P)$ and $\mathbf{\Psi }(Q)$ are easily computable and can be viewed as points of $\mathbb{P}_{\C }^7$.  The matrix $\mathbf{\Psi }(P)$ is hyperbolic-like, however, the matrix $\mathbf{\Psi }(Q)$ is not. The maximal eigenvalue (in absolute value) of it has multiplicity two, and so is the minimal eigenvalue.

  \medskip 
  
  We consider the following conditions for $\mathbf{\Psi }(P)$ and $\mathbf{\Psi }Q)$.

  \medskip 
  
  (1) $\mathcal{F}_{\mathbf{\Psi }(P)} \cap \mathcal{F}_{\mathbf{\Psi }(Q)} = \emptyset $; 

  \medskip 

  (2) $\mathcal{F}_{\mathbf{\Psi }(P) }\cap \Pi _{\mathbf{\Psi }(Q)}= \emptyset = \mathcal{F}_{\mathbf{\Phi }(Q)} \cap \Pi _{\mathbf{\Psi }(P)} $;

  \medskip 
  
  (3) for all sufficiently big $M, N\in \N $, $[\mathbf{\Psi }(Q)^M  \mathbf{\Psi }(P)^N \mathbf{\Psi }(Q)^{-M}, \mathbf{\Psi }(P)^N]\neq 1$ and $[\mathbf{\Psi }(P)^M  \mathbf{\Psi }(Q)^N \mathbf{\Psi }(P)^{-M}, \mathbf{\Psi }(Q)^N]\neq 1$. 

    \medskip 

    We cannot claim conditions (1)-(2), but we will find an easy substitute below, more precisely, we have these conditions satisfied for matrices $\sigma _2(P)$ and $\sigma _2(Q)$. In other words, the matrices $\sigma _2(P)$ and $\sigma _2(Q)$ are hyperbolic without a common eigenvector. Then by Proposition \ref{thm:lastfactor}, for sufficiently big $N$, for all $m, n\geq N$ the matrices $P^m$ and $Q^n$ generate a non-Abelian free group. 
    
    \medskip 
    
    Notice that, since $\mathbf{\Psi }$ is a monomorphism, condition (3) is implied by the conditions $$[(\sigma _2Q)^M(\sigma _2P)^N(\sigma _2Q)^{-M},(\sigma _2P)^N]\neq 1 \ \mathrm{and} \  [(\sigma _2P)^M(\sigma _2Q)^N(\sigma _2P)^{-M},(\sigma _2Q)^N]\neq 1$$ for all non-zero integer $M$. Thus, it is straightforward (a direct computation) to satisfy conditions (1)-(3). For condition (3), alternatively, for sufficiently large $M$ and $N$, it immediately follows from condition (1) and Proposition \ref{prop:noncommuting}.

    \medskip

    Now we are ready for a quick finishing argument. We let $f = (P, \sigma _1(P)), g = (Q, \sigma _1(Q))$ and $\Gamma := \Gamma _N = \langle f^N, g^N\rangle $.

    \medskip 

    Let also $\pi _0:SL(2,\R)\times SL(2, \C)\to SL(2, \R), \pi _1:SL(2,\R)\times SL(2, \C)\to SL(2, \C)$ be the projections onto the first and second coordinates. Then $\pi _0(\Gamma ) = \langle P^N, Q^N \rangle $ and $ \pi _1(\Gamma ) = \langle \sigma _1(P)^N, \sigma _1(Q)^N \rangle $. 

    \medskip

        Notice that the matrix $P$ is elliptic but not a torsion. We claim that for sufficiently large $N$, the subgroup $\langle P^N, Q^N \rangle $ is non-Abelian free. Indeed, if not, then there is a non-trivial relation between $\mathbf{\Psi }(P)^N$ and $\mathbf{\Psi }(Q)^N$. Then the same relation holds for $\sigma _2(P)^N$ and $\sigma _2(Q)^N$. The latter is a pair of hyperbolic $2\times 2$ matrices, and direct computation shows that conditions (1) and (2) hold for these matrices (in size two, these two conditions become equivalent).

\medskip 

         Thus, we established that $\Gamma $ is non-Abelian free. In the next section, we will verify that $\Gamma $ is discrete. Then, it remains to show that the projections of $\Gamma $ to first and second coordinates are both dense. Since $P$ is an elliptic non-torsion element, the closure $\overline{\langle P^N\rangle }$ is isomorphic to $\mathbb{S}^1$. By condition (2), the closure $\overline{\pi _0(\Gamma )} = \overline{\langle P^N, Q^N\rangle }$, as a Lie subgroup, contains infinitely many copies of $\mathbb{S}^1$, hence it is at least two-dimensional (let us also recall a classical fact, due to E.Cartan, that a closed subgroup of a Lie group is a Lie subgroup). On the other hand, since two-dimensional connected Lie groups are solvable (hence they do not contain a copy of $\mathbb{F}_2$), the closure $\overline{\langle P^N, Q^N\rangle }$, as a Lie subgroup, must be 3-dimensional. Hence $\overline{\langle P^N, Q^N\rangle } = SL(2, \R)$. 

         \medskip 
         
         Similarly, we can claim that $\overline{\langle \sigma _1(Q)^N, \sigma _1(P)^N\rangle } = SL(2, \C)$. For this, in addition, we observe that $\overline{\langle \sigma _1(Q)^N, \sigma _1(P)^N\rangle }$ is not compact thus it is not contained in any conjugate of $SU(2)$. By looking at the trace, we also conclude that $\overline{\langle \sigma _1(Q)^N, \sigma _1(P)^N\rangle }$  is not conjugate to a subgroup of $SL(2,\R)$. Then the closure $\overline{\langle \sigma _1(Q)^N, \sigma _1(P)^N\rangle }$ is a non-compact Lie subgroup of $SL(2,\C )$ of dimension at least three other than a conjugate copy of $SL(2,\R)$.   Hence $\overline{\langle \sigma _1(Q)^N, \sigma _1(P)^N\rangle } = SL(2, \C)$. 
    
 \vspace{1cm}

\section{Verifying discreteness of $\Gamma $}

 Our group $\Gamma $ is generated by $\langle P^N, \sigma _1(P)^N\rangle $ and $\langle Q^N, \sigma _1(Q)^N\rangle $. We will consider extended groups $$\Gamma _1 := \langle (P, \sigma _1(P), \sigma _3(P) ) , (Q, \sigma _1(Q), \sigma _3 (Q) ) \rangle $$ \ $$ \Gamma _2:= \langle (P, \sigma _1(P), \sigma _2(P), \sigma _3(P)), (Q, \sigma _1(Q), \sigma _2(Q), \sigma _3(Q))\rangle $$ and $$\Gamma _3 := \langle (P, \sigma _1(P), \sigma _2(P)), (Q, \sigma _1(Q), \sigma_2(Q)) \rangle .$$

   \medskip 

   The generators of $\Gamma _1$ can be presented as the triples $$f_1 := (\begin{bmatrix} 5+ \beta ^2 - 3\beta -2\beta ^3 & 1 \\
	-1 & 0 
        \end{bmatrix},  \begin{bmatrix} 5- \beta ^2 - (3\beta -2\beta ^3){\bf i}& 1 \\
	-1 & 0 
        \end{bmatrix}, \begin{bmatrix} 5- \beta ^2 + (3\beta -2\beta ^3){\bf i} & 1 \\
	-1 & 0 
        \end{bmatrix} )$$ and $$g_1 := (\begin{bmatrix} 3+ 2\beta ^2  & 1 \\
	-1 & 0 
        \end{bmatrix},  \begin{bmatrix} 3- 2\beta ^2  & 1 \\
	-1 & 0 
        \end{bmatrix}, \begin{bmatrix} 3 -2\beta ^2  & 1 \\
	-1 & 0 
        \end{bmatrix}).$$

\medskip 

The group $\Gamma _3$ is discrete in $SL(2,\R)\times SL(2,\C) \times SL(2,\R)$, and the group $\Gamma _2$ is discrete in $SL(2,\R)\times SL(2,\C) \times SL(2,\R)\times SL(2,\C)$. On the other hand, we do not know yet if the group $\Gamma _1$ is necessarily discrete in $SL(2,\R)\times SL(2,\C) \times SL(2,\C)$.  The generators of $\Gamma _2$ can be presented as quadruples $$\footnotesize{ (\begin{bmatrix} 5+ \beta ^2 - 3\beta -2\beta ^3 & 1 \\
	-1 & 0 
        \end{bmatrix},  \begin{bmatrix} 5- \beta ^2 - (3\beta -2\beta ^3){\bf i}& 1 \\
	-1 & 0 
        \end{bmatrix}, \begin{bmatrix} 5+ \beta ^2 +3\beta +2\beta ^3 & 1 \\
	-1 & 0 
        \end{bmatrix}, \begin{bmatrix} 5+ \beta ^2 + (3\beta -2\beta ^3){\bf i} & 1 \\
	-1 & 0 
        \end{bmatrix} )}$$

and $$(\begin{bmatrix} 3+ 2\beta ^2  & 1 \\
	-1 & 0 
        \end{bmatrix},  \begin{bmatrix} 3- 2\beta ^2  & 1 \\
	-1 & 0 
        \end{bmatrix}, \begin{bmatrix} 3+ 2\beta ^2  & 1 \\
	-1 & 0 
        \end{bmatrix}, \begin{bmatrix} 3 -2\beta ^2  & 1 \\
	-1 & 0 
        \end{bmatrix}).$$ Let $f_2, g_2$ be these generators respectively.

\medskip

   By abuse of notation, we also let $\pi _k:\Gamma _2\to SL(2,\C), k\in \{0,1,2,3\}$  and $\pi _l:\Gamma _3\to SL(2,\C), l\in \{0,1,2\}$ be the projections to the $k$-th factor and $l$-th factor respectively. Notice that for $k\in \{0,2\}$ and $l\in \{0,2\}$ the image of the projection lies in $SL(2,\R)$. For $i\in \{1,2,3\}$, we will also let $\Gamma _{i,N}$ denote the group generated by $f_i^N$ and $g_i^N$. 

   \medskip 
   
   For any two matrices $A = (a_{ij})_{1\leq i,j\leq n}, B = (b_{ij})_{1\leq i,j\leq n}$ in $\mathrm{Mat}(n,\C )$, we also let $d(A,B) := \displaystyle \mathop{\max }_{1\leq i, j\leq n}|a_{ij}-b_{ij}|$. This metric can be extended naturally to any product $\mathrm{Mat}(n_1, \C )\oplus \dots \oplus \mathrm{Mat}(n_k, \C)$ taking the supremum of distances in each coordinate. For all $A\in \mathrm{Mat}(n,\C )$, we also let $||A|| = d(A, I) = \displaystyle \mathop{\max }_{1\leq i, j\leq n}|a_{ij}-\delta _i^j|$ and $||A||_0 = \displaystyle \mathop{\min }_{1\leq i, j\leq n}|a_{ij}-\delta _i^j|$. Also, for $1\leq i, j\leq n$, we let $Ent_{ij}(A)$ denote the $(ij)$-entry of $A$, and we let $\mathrm{Re}(A)$ and $\mathrm{Im}(A)$ denote the real part and the imaginary part of the matrix $A$ respectively.   

 \medskip 

  In the ring $\Z[\beta ]$, we introduce the quantity $$N(x) = |\sigma _0(x)\sigma _1(x)\sigma _2(x)\sigma _3(x)|.$$ It would be useful to recall that $\sigma _0(x) = x$ and $\sigma _1(x) = \overline{\sigma _3(x)}$ for any $x\in \Z[\beta ]$. If $x = a + m\beta + p\beta ^2 + e\beta ^3, a,m,p,e\in \Z$, then one can compute that $$N(x) = |(a^2+2p^2-4me)^2-2(2ap-m^2-2e^2)^2|$$ thus $N(x) \geq 1$ unless $x=0$. \footnote{The analog of the quantity $N(x)$  can be defined in other Galois rings as well. In the ring $\Z [\sqrt{2}]$, we can define it as $N(x) = |(m+n\sqrt{2})(m-n\sqrt{2})|$ for $x = m+n\sqrt{2}, m.n\in \Z$. Again, we observe that $N(x) = |m^2-2n^2| \geq 1$ unless $x=0$.} An element $x = a + m\beta + p\beta ^2 + e\beta ^3$ of $\Q[\beta ]$ will be called {\em positive}, if $a,m,p,e\geq 0$ and $x\neq 0$; $x$ is called {\em negative} if $-x$ is positive; and $x$ is called a {\em signed element} if it is either positive or negative. For a positive $x$, we also let $$\Delta (x) = x - \gamma (x), \ \mathrm{and} \ \Delta _i(x) = x_i - \gamma _i(x), 1\leq i\leq 2$$ where $$\gamma (x) = 4\min \{a, m\beta , p\beta ^2 , e\beta ^3\}, \gamma _1(x) = 2\min \{a, p\beta ^2\}, \gamma _2(x) = 2\min \{m\beta , e\beta ^3\},$$ and  $x_1 = a+p\beta ^2, x_2= m\beta +e\beta ^3.$ The quantities $\Delta , \Delta _1, \Delta _2$ and $\gamma , \gamma _1, \gamma _2$ can also be extended to negative elements as well by letting, for all $1\leq i\leq 2$,  $$\Delta (x) = -\Delta (-x), \Delta _i(x) = -\Delta _i(-x)$$ and $$\gamma (x) = -\gamma (-x), \gamma _i(x) = -\gamma _i(-x)$$ when $x$ is negative. Thus, the equation $$x = \Delta (x) + \gamma (x) = \Delta _1(x) + \gamma _1(x) + \Delta _2(x) + \gamma _2(x)$$ also holds for negative $x$. 

  \medskip 

  If $x, y\in \Q[\beta ]$ are signed elements, then $xy$ is also signed, moreover, $xy$ is positive if and only if $x$ and $y$ have the same signs. Also, for a signed $x\in \Z[\beta ]$, $\gamma (x)$ is also signed (with the same sign), but $\Delta (x)$ is not (because of the irrationality of $\beta, \beta ^2$ and $\beta ^3$). However, for any $\delta > 0$, there exists a signed $x' = a' + m'\beta + p'\beta ^2 + e'\beta ^3 \in \Q[\beta ]$ such that $|x'-x| < \delta $, $\mathrm{diam} \{a', m'\beta , p'\beta ^2, e'\beta ^3, \frac{\gamma (x)}{4}\} < \delta $ and $x', \Delta (x')$ are signed elements of $\Q[\beta ]$ with the same sign as of $x$; hence we also have $|\gamma (x)-\gamma (x')| < 4\delta $.  Then for all signed $x, y\in \Z[\beta ]$, from  $xy = $ \ $$(\gamma (x)+\Delta (x))(\gamma (y)+\Delta (y)) = (\gamma (x)\gamma (y) + \gamma (x)\Delta (y) + \gamma (y)\Delta (x)) + \Delta (x)\Delta (y)$$ we obtain $|\Delta (xy)|\leq |\Delta (x)\Delta (y)| \ (1)$. 
  
  \medskip 
  
  Inequality (1) as an indication of some coherence in the behavior of $\Delta $; its analog can also be stated for $\Delta _1$ and $\Delta _2$. However, we will need to deal with products $xy$ where only one of the factors is signed. To deal with this case, we will make the following observations. For positive real numbers $\epsilon $ and $c$, let $$S_{\epsilon , c} = \{x\in \Z[\beta ] \ | \ 0 < |x| < \epsilon , |\sigma _1(x)| < c \}.$$ We make a useful observation that for a fixed $c> 0$, if $\epsilon > 0$ is sufficiently small, then for all $x\in S_{\epsilon , c}$, $\sigma _2(x)$ is a signed element, $$|\Delta _1(\sigma _2(x))| = |\mathrm{Re}\sigma _1(x)|, |\Delta _2(\sigma _2(x)) = \mathrm{Im}\sigma _1(x)|$$ and by the inequality $N(x)\geq 1$, we have $|\sigma _2(x)| \geq \frac{1}{\epsilon c^2}$ thus, $\displaystyle \mathop{\lim }_n|\sigma _2(x)| = \infty $ as $\epsilon \to 0$ and $c$ is fixed. But we make a stronger observation that even for a broader set $S_{c,c}$ with a fixed $c >0$ if $(x_n)$ is a sequence of distinct elements in it, then for sufficiently large $n$, the element $\sigma _2(x_n)$ is signed and $\displaystyle \mathop{\lim }_n|\sigma _2(x_n)| = \infty $.  However, we make another very useful observation that for all fixed $c > 0$, if $(x_n)$ is a sequence of distinct elements in  $$S_{c,c}' := \{x\in \Z[\beta ] \ | \ 0 < |x| < c , \mathrm{Re}(\sigma _1(x)) \ \mathrm{or} \ \mathrm{Im}(\sigma _1(x)) \ \mathrm{is \ not \  signed} \},$$ then $\displaystyle \mathop{\lim }_n|\sigma _2(x_n)| = \infty $ provided that both $(\mathrm{Re}(\sigma _1(x_n)))$ and $(\mathrm{Im}(\sigma _1(x_n)))$ are sequences of distinct elements, i.e., if $m\neq n$, then $\mathrm{Re}(\sigma _1(x_m))\neq \mathrm{Re}(\sigma _1(x_n))$ and $\mathrm{Im}(\sigma _1(x_m))\neq \mathrm{Im}(\sigma _1(x_n))$).

  \medskip

 Let us now assume that $\Gamma $ is not discrete in $SL(2,\R)\times SL(2,\C)$. Then for all $\epsilon > 0$ there exists a non-identity matrix $$A = \begin{bmatrix} a + m\beta  + p\beta ^2 + e\beta ^3    & b + n\beta  + q\beta ^2 + f\beta ^3  \\
	c + k\beta  + r\beta ^2 + g\beta ^3  & d + l\beta  + s\beta ^2 + h\beta ^3  
        \end{bmatrix} $$ in $\langle P^N, Q^N\rangle $ such that $d(A, I) < \epsilon  $ and $d(\sigma _1(A), I) < \epsilon $ where $a, m, \dots , s, h$ are integers. Then $d(\sigma _3(A), I) < \epsilon $  as well, and we can claim that $\Gamma _{1, N} := \langle f_1^N, g_1^N\rangle $ is also non-discrete. On the other hand, the inequality $d(\sigma _1(A), I) < \epsilon $  implies that $$\max \{|a-1-p\beta ^2| , |m-e\beta ^2|,  |b-q\beta ^2|, |n-f\beta ^2| ,  |c-r\beta ^2|, |k-g\beta ^2|, |d-1-s\beta ^2|, |l-h\beta ^2| \} < \epsilon .$$  Then the inequality $d(A,I) < \epsilon $ implies that the numbers in the quadruple $(a-1, -m\beta , p\beta ^2, -e\beta ^3)$ are at most $4\epsilon $ apart.  Similarly, in each of the quadruples $$(b, -n\beta , q\beta ^2, -f\beta ^3),   (c, -k\beta , r\beta ^2, -g\beta ^3),  (d-1, -l\beta , s\beta ^2, -h\beta ^3),$$ any two coordinates are at most $4\epsilon $ apart. In addition, we also have  $(\Psi A) = $ \ $$\begin{bmatrix} 1 - e\beta ^3 &  2e & - 2e\beta  & 2e\beta ^2  & - f\beta ^3 &  2f & -2f\beta & 2f\beta ^2  \\
e\beta ^2 & 1 - e\beta ^3 & 2e & -2e\beta & f\beta ^2 & - f\beta ^3 & 2f & -2f\beta \\ -e\beta & e\beta ^2 & 1-e\beta ^3 & 2e  & -f\beta & f\beta ^2 & -f\beta ^3 & 2f \\
e & -e\beta & e\beta ^2 & 1-e\beta ^3 & f & -f\beta & f\beta ^2 & -f\beta ^3 \\ 

- g\beta ^3 &  2g & - 2g\beta  & 2g\beta ^2  & 1- h\beta ^3 &  2h & -2h\beta & 2h\beta ^2  \\
g\beta ^2 & - g\beta ^3 & 2g & -2g\beta & h\beta ^2 & 1- h\beta ^3 & 2h & -2h\beta \\ -g\beta & g\beta ^2 & -g\beta ^3 & 2g  & -h\beta & h\beta ^2 & 1-h\beta ^3 & 2h \\
g & -g\beta & g\beta ^2 & -g\beta ^3 & h & -h\beta & h\beta ^2 & 1-h\beta ^3  
        \end{bmatrix}$$ For arbitrary $\epsilon > 0$ and $M>0$, we can also assume that $||A|| < \epsilon $ and  $\max \{|e|,|f|,|g|,|h|\} > M$.  The latter implies that either $\max \{|e|,|g|\} > M$ or $\max \{|f|,|h|\} > M$; without loss of generality we will assume that $\max \{|e|,|g|\} > M$. Then, since $|\det(\Psi (A))| = 1$, we can also assume that $\max \{|f|,|h|\} > M$. Thus, for all entries $x$ of the matrix $A$, we can claim that $|\sigma _2(x)| > M$. Again, we can make a stronger observation that if $(A_k)$ is a sequence of matrices with entries in $S_{c,c}'$ for some fixed $c>0$ (rather than $S_{\epsilon , \epsilon }$ when $\epsilon \to 0$), then we can still claim that if any two entries of $A_k$ and $A_l$ are distinct for $k\neq l$, then for sufficiently large $n$, the inequality $|\sigma _2(x) > M$ holds for any entry of $x$ of $A_n$.   

 \medskip

 Let $u_1= [1, \lambda _1], w_1 = [1, \lambda _2]$ be the eigenvectors of $\sigma _2(P)$ and $u_2 = [1, \lambda _3], w_2= [1, \lambda _4]$ be the eigenvectors of $\sigma _2(Q)$. We will make use of the fact that $|\lambda _i|\neq 1, 1\leq i\leq 4$. If $v$ is a fixed vector which is not collinear with these four vectors (i.e. $[v]\notin \{[u_1], [w_1], [u_2], [w_2]\}$ in $\C P^1$), for sufficiently large $N$, $[Aw]$ will be close to one of the points $[u_1], [w_1], [u_2], [w_2]$ in $\C P^1$. Then for all $D > 0$, since $M > 0$ can be chosen sufficiently large (and $\epsilon $ sufficiently small), without loss of generality we may assume that
$$\mathrm{dist} ([e:g], [1:\lambda_i])<D \ \mathrm{and} \ \mathrm{dist} ([f:h], [1:\lambda_i])<D \ (2)$$ for some $1\leq i\leq 4$.

\medskip

 Let $$\sigma _2(A) = \begin{bmatrix} \zeta  & \eta \\ \mu & \nu \end{bmatrix}, \sigma _1(A) = \begin{bmatrix} \zeta ' & \eta '\\ \mu '& \nu '\end{bmatrix}, \ \mathrm{and} \ A = \begin{bmatrix} \zeta '' & \eta ''\\ \mu ''& \nu ''\end{bmatrix} $$  \footnote{Since $\sigma _2$ is an involution, we have $\zeta '' = \sigma _2(\zeta ), \eta '' = \sigma _2(\eta ), \mu '' = \sigma _2(\mu )$ and $\nu '' = \sigma _2(\nu )$.}

 \medskip 
 
  Since $\Gamma $ is free, both factors $\pi _0(\Gamma )$ and $\pi _1(\Gamma )$ are also free, therefore these factors are torsion-free. Since $\Gamma $ is not discrete, its closure will be a Lie subgroup of $SL(2,\R)\times SL(2,\C)$ containing a copy of $SL(2,\R)$ in the first coordinate and a copy of $SL(2,\C)$ in the second coordinate. Considering Lie subgroups of $SL(2,\R)\times SL(2,\C)$, similar to the argument at the end of Section 3, we conclude that the closure of $\Gamma $ will be equal to $SL(2,\R)\times SL(2,\C)$. 

  \medskip 

  The group $\pi _2(\Gamma _{3,N})$ consists of the elements represented by the reduced words in the alphabet $\{\sigma _2(P)^N, \sigma _2(Q)^N\}$. Let $\Omega _N$ denote the set of hyperbolic matrices in $\pi _2(\Gamma _{3,N})$ represented by the reduced words $W$ such that if $W$ starts with $\sigma _2(P)^{\pm N}$, then it ends with $\sigma _2(Q)^{\pm N}$, and if $W$ starts with $\sigma _2(Q)^{\pm N}$, then it ends with $\sigma _2(P)^{\pm N}$. An important property of the set $\Omega _N$ for us will be the fact that for sufficiently large $N$, the two eigenvectors of the matrices in $\Omega _N$ stay uniformly away from each other in the standard metric of $\mathbb{CP}^{1}$, i.e. there exists $D > 0$ such that for all $A\in \Omega _N$, if $v_1, v_2$ are the two eigenvectors of $A$, then $\mathrm{dist}(v_1, v_2) > D$. 

  \medskip 

   By inequality (2), we may assume that $$|\frac{\zeta }{\mu }|, |\frac{\eta }{\nu }|\in \displaystyle \mathop{\cup }_{i=1}^{4}(\lambda _i-D, \lambda _i + D).$$

   Motivated by this observation, for all $\kappa > 0, D > 0$, we introduce a {\em cone} $$\mathcal{C}_{\kappa, D} = \{\begin{bmatrix} \zeta  & \eta \\ \mu & \nu \end{bmatrix}\in \pi _2(\Gamma _{3,N}) : |\frac{\zeta }{\mu }|, |\frac{\eta }{\nu }|\in (\kappa -D, \kappa  + D)\}.$$ 
   
   We make another observation that given any $\delta > 0$, for sufficiently large $N$, the projection of the group $\Gamma = \Gamma _N$ (more precisely, the set $\pi _2(\Gamma _{3,N})\backslash \{1\}$) will lie in the union of the cones $\mathcal{C}_{\lambda _1, \delta }, \mathcal{C}_{\lambda _1^{-1}, \delta }, \mathcal{C}_{\lambda _2, \delta }, \mathcal{C}_{\lambda _2^{-1}, \delta }$. 

 \medskip 

  It will be useful to broaden the notion of a cone to include infinite values of $\kappa $; a diagonal matrix (in particular, a matrix $\begin{bmatrix} \lambda ^n & 0 \\ 0 & \lambda ^{-n}  \end{bmatrix}$ for any integer $n$) can be viewed as an element of $\mathcal{C}_{\infty , D}$ for any $D>0$.). 

  \medskip 

  Notice that if $B = \begin{bmatrix} \lambda ^n & 0 \\ 0 & \lambda ^{-n}  \end{bmatrix}, \lambda \in \C, |\lambda | > 1$ and $C = \begin{bmatrix} a & b \\ c & d  \end{bmatrix}\in SL_2(\C )$, then $CBC^{-1} = \begin{bmatrix} ad\lambda ^n-bc\lambda ^{-n} & ab(\lambda ^{-n}-\lambda ^n) \\ cd(\lambda ^n-\lambda ^{-n}) & ad\lambda ^{-n}-bc\lambda ^n  \end{bmatrix}$. This implies that if no entries of $C$ are zero, then for all fixed $\delta > 0$, if $n$ is sufficiently large, we have $CBC^{-1}\in \mathcal{C}_{a/c, \delta }$. More generally, if $B = \begin{bmatrix} \lambda  & 0 \\ 0 & \lambda   \end{bmatrix}, \lambda \in \C, |\lambda | > 1$, then, given $\delta > 0$, for sufficiently large  $\lambda $, $CBC^{-1}\in \mathcal{C}_{a/c, \delta }$ and $CB^{-1}C^{-1}\in \mathcal{C}_{b/d, \delta }$. On the other hand, for all $z\in \C\backslash \{0\}$, $$\begin{bmatrix} z & 0 \\ 0 & z^{-1} \end{bmatrix}B\begin{bmatrix} z^{-1} & 0 \\ 0 & z \end{bmatrix} = B \ \mathrm{and} \ C\begin{bmatrix} z & 0 \\ 0 & z^{-1} \end{bmatrix} = \begin{bmatrix} za & z^{-1}b \\ zc & z^{-1}d \end{bmatrix}.$$ We can choose $z\in \C$ such that the sum $$S(z):=\frac{|za|}{|z^{-1}d|}+\frac{|z^{-1}d|}{|za|} + \frac{|z^{-1}b|}{|zc|}+\frac{|zc|}{|z^{-1}b|}$$ is minimal; and if $\min S$ is such minimal value, by density of $\Z[\beta ]$ in $\R $, for all $\delta > 0$, we can also choose $z'\in \Z[\beta ]$ such that $S(z') < \min S + \delta $.
  
  Since the set $$ \{\begin{bmatrix} a & b \\ c & d  \end{bmatrix}\in SL_2(\C ) : |\frac{a}{c}|+|\frac{c}{a}|+|\frac{b}{d}|+|\frac{d}{b}|+|\frac{a}{d}|+|\frac{d}{a}|\leq M\}$$  
   is compact for all $M > 0$, by density of $\Gamma $, we can choose sequences $C_k\in SL(2,\R), k\geq 1$ and $A_k\in \pi _2(\Gamma _{3,N}), k\geq 1$ and a matrix $B\in \pi _2(\Gamma _{3,N})$ such that the following conditions hold: 

   \medskip 

   1) $A_k\in \Omega _N, k\geq 1$ and $\sigma _2(A_k), k\geq 1$  are elliptic matrices \footnote{Here and in the sequel, maybe it is more natural to write $\sigma _2^{-1}$ instead of $\sigma _2$, but notice that $\sigma _2$ is an involution.}; 

   \medskip 
   
   2) $\sigma _2(A_k)\to A$ as $k\to \infty $ for some elliptic matrix $A\in SL(2,\R )$ and the sequence $(\sigma _1^{-1}(A_k))$ remains in a compact subset of $SL(2,\R)$;
   
   \medskip 
   
   3) $\displaystyle \mathop {\inf }\{||\mathrm{Re}(\sigma _1^{-1}A_k)||_0 :k\geq 1\} > 0$ and $\displaystyle \mathop {\lim }_{k\geq 1} ||\mathrm{Im}(\sigma _1^{-1}A_k)|| = 0$; moreover, each of the sequences $(\frac{Ent_{11}(\mathrm{Re}(\sigma _1^{-1}A_k))}{Ent_{12}(\mathrm{Re}(\sigma _1^{-1}A_k))})$ and $(\frac{Ent_{21}(\mathrm{Re}(\sigma _1^{-1}A_k))}{Ent_{22}(\mathrm{Re}(\sigma _1^{-1}A_k))})$ is dense in some interval of positive length (not necessarily the same interval); 

\medskip 

   4) if $x, y$ are entries in $A_k$ and $A_l$, respectively, for some $k\neq l$, then $x\neq y$; 

\medskip    

  5) the sequence $(C_k)$ remains in a compact subset of $SL(2,\R)$, and for all $k\geq 1$, $C_k^{-1}A_kC_k$ is a diagonal matrix with real entries $\theta _k, \theta _k^{-1}$ and $\theta _k > 1 > \theta _k^{-1}$ (so,  $\theta _k, \theta _k^{-1}$ are the eigenvalues of $A_k$); 
  
  \medskip 

  6)  $\sigma _2(B)$ is an elliptic matrix, and there exists a matrix $C'\in SL(2,\R)$ such that $(C')^{-1}BC'$ is a diagonal matrix with real entries $\theta, \theta ^{-1}$ and $\theta > 1 > \theta ^{-1}$; 

  \medskip 

  7) for all $k\geq 1$ at least one of the elements of the set $$\{Ent_{ij}\mathrm{Re}(\sigma _1^{-1}(A_kB^{-q_k})) : i,j\in \{1,2\} \}\cup \{Ent_{ij}\mathrm{Im}(\sigma _1^{-1}(A_kB^{-q_k})) : i,j\in \{1,2\}\}$$ is not a signed element, where $q_k = [\log \frac {\theta _k}{\theta }]$.

  \medskip  

   Notice that by condition 2) and 4), $Tr(A_k) \to \infty $ as $k\to \infty $, therefore $A_k$ is a hyperbolic matrix for sufficiently large $k$ with its larger eigenvalue tending to infinity. This, together with condition 1) (since $A_k\in \Omega _N, k\geq 1$, the eigenvectors of $A_k$ stay uniformly away from each other) allows us to choose $(C_k)$ from a compact subset to guarantee condition 5). By passing to a subsequence if necessary, we can also assume that $C_k\to C$ for some $C\in SL(2,\R)$.

   \medskip 
   
   Regarding condition 7), let us point out that for the matrix $$P =  \begin{bmatrix} 
	5+ \beta ^2 - 3\beta -2\beta ^3 & 1 \\
	-1 & 0 
 \end{bmatrix},$$ for sufficiently large $k$, at least one of the entries of $\mathrm{Re} \sigma _1(P)^k$ or $\mathrm{Im} \sigma _1(P)^k$ is not a signed element. This is because $P$ itself is elliptic (so $||P^k||, k\geq 1$ is a bounded sequence) and $\displaystyle \mathop{\lim }_{k\to \infty }\frac{||\sigma _2(P)^k||}{||\sigma _1(P)^k||} = \infty $. However, by a straightforward induction, it is easy to see that for all $k\geq 1$, all four entries of $\sigma _2(P)^k$ are signed elements.\footnote{Notice that $\sigma _2(P)^k = \begin{bmatrix}  u_{k+2} & u_{k+1} \\ -u_{k+1} & u_{k}  \end{bmatrix}$ for all $k\geq 1$, with $u_1=0, u_2 = 1, u_{k+2}=\omega u_{k+1}-u_k$ and $\omega = 5+\beta ^2+3\beta +2\beta ^3$.} Then every entry in $\mathrm{Re} \sigma _1(P)^k$ and in $\mathrm{Im} \sigma _1(P)^k$ is not a signed element.

\medskip

   Notice that $C_k\begin{bmatrix} \theta _k & 0 \\ 0 & \theta _k^{-1}  \end{bmatrix}C_k^{-1}, k\geq 1$ and $C'\begin{bmatrix} \theta  & 0 \\ 0 & \theta ^{-1}  \end{bmatrix}(C')^{-1}$ are matrices in $SL(2,\Z[\beta ])$.
   
 \medskip

Then we can choose a sequence $D_k\in SL(2, \Z[\beta )$) such that $[D_k,\sigma _2(P)]=1, k\geq 1$, and the sequence $(D_k)$ does not leave a compact set $K$ in $SL(2,\R)$; moreover, $\displaystyle \mathop{\lim }_{k\to \infty }E_k||A_k||^2 = 0$ where $E_k = Ent_{12}(C_k^{-1}D_kC'), k\geq 1$ (so, in particular, the sequence $E_k$ converges to zero as $k\to \infty $).  To arrange the commuting relations, it suffices to take $D_k$ in the form $\begin{bmatrix} \omega y+t & y \\ -y & t]  \end{bmatrix}$ where $y, t\in \Z[\beta ], \omega = 5 + 3\beta + \beta ^2 + 2 \beta ^3$ and $t^2+y^2+\omega yt = 1$ (i.e. the determinant equals 1). Viewing the latter as a quadratic equation in the variable $t$ over the ring $\Z[\beta ]$ leads (by arranging the discriminant as a perfect square in this ring) to a Pell-type equation $s^2-(\omega ^2-4)y^2=4$.  This equation is a rank two Pell-type equation and invoking logarithmic independence one can see that it has a solution for infinitely many $s$ from a bounded interval (cf. \cite{Sch}), thus we can also arrange the condition $\displaystyle \mathop{\lim }_{k\to \infty }E_k||A_k||^2 = 0$.

\medskip 

Notice that, by conditions 2) and 6), for any sequence $(p_k)$ of integers, the sequence $$\sigma _2((C_k\begin{bmatrix} \theta _k & 0 \\ 0 & \theta _k^{-1}  \end{bmatrix}C_k^{-1})(C'\begin{bmatrix} \theta  & 0 \\ 0 & \theta ^{-1}  \end{bmatrix}^{-p_k}(C')^{-1})) = \sigma _2(A_kB^{-p_k})$$ stays in a compact subset of $SL(2,\R)$. 

 \medskip 
 
 On the other hand, for the choice $p_k = q_k, k\geq 1$ and $B = \sigma _2(P)$, and using conditions 3) and 7), we can achieve that at least one of the entries of either real or imaginary part of the matrix $$\sigma _1^{-1}((C_k\begin{bmatrix} \theta _k & 0 \\ 0 & \theta _k^{-1}  \end{bmatrix}C_k^{-1})(C'\begin{bmatrix} \theta  & 0 \\ 0 & \theta ^{-1}  \end{bmatrix}^{-p_k}(C')^{-1})) = \sigma _1^{-1}(A_kB^{-p_k})$$ will not be signed. 

  \medskip
  
  Therefore, since $$\sigma _2((C_k\begin{bmatrix} \theta _k & 0 \\ 0 & \theta _k^{-1}  \end{bmatrix}C_k^{-1})(C'\begin{bmatrix} \theta  & 0 \\ 0 & \theta ^{-1}  \end{bmatrix}^{-p_k}(C')^{-1}))$$ remains in a compact, the sequence $$C_k\begin{bmatrix} \theta _k & 0 \\ 0 & \theta _k^{-1}  \end{bmatrix}C_k^{-1}C'\begin{bmatrix} \theta  & 0 \\ 0 & \theta ^{-1}  \end{bmatrix}^{-q_k}(C')^{-1} = A_kB^{-q_k}$$ escapes every compact subset. Then, since $(D_k)$ remains in a compact, and $D_k$ commutes with $B$ for all $k\geq 1$, the sequence $$C_k\begin{bmatrix} \theta _k & 0 \\ 0 & \theta _k^{-1}  \end{bmatrix}C_k^{-1}D_kC'\begin{bmatrix} \theta  & 0 \\ 0 & \theta ^{-1}  \end{bmatrix}^{-q_k}(C')^{-1} = C_k\begin{bmatrix} \theta _k & 0 \\ 0 & \theta _k^{-1}  \end{bmatrix}C_k^{-1}C'\begin{bmatrix} \theta  & 0 \\ 0 & \theta ^{-1}  \end{bmatrix}^{-q_k}(C')^{-1}D_k$$ also escapes from every compact subset.
  
  \medskip 
  
  However, notice that the matrix $C_k^{-1}D_kC'$ is of the form $\begin{bmatrix} a_k  & E_k \\ b_k & c_k  \end{bmatrix}$ where the sequences $(a_k), (b_k), (c_k)$ are bounded and $\lim E_k||A_k||^2 = 0$. But  $$\begin{bmatrix} \theta _k & 0 \\ 0 & \theta _k^{-1}  \end{bmatrix}C_k^{-1}D_kC'\begin{bmatrix} \theta  & 0 \\ 0 & \theta ^{-1}  \end{bmatrix}^{-p_k} = \begin{bmatrix} \theta _ka_k\theta ^{-p_k}  & E_k\theta _k\theta ^{p_k} \\ \theta _k^{-1}c_k\theta ^{-p_k} & \theta _k ^{-1}d_k\theta ^{p_k} \end{bmatrix}$$ therefore, for a suitable choice of $p_k$ (namely, for $p_k = [\log \frac{\theta _k}{\theta }]$) all entries of the latter matrix remain bounded. Then all entries in the matrix $$C_k\begin{bmatrix} \theta _k & 0 \\ 0 & \theta _k^{-1}  \end{bmatrix}C_k^{-1}D_kC'\begin{bmatrix} \theta  & 0 \\ 0 & \theta ^{-1}  \end{bmatrix}^{-q_k}(C')^{-1}D_k^{-1} = A_kB^{-q_k}$$ also remain bounded. This is a contradiction. 

   \medskip 

 Thus, the subgroup $\Gamma $ is discrete in $SL(2,\R)\times SL(2,\C)$.

\bigskip

 \section{Proof of Theorem \ref{thm:main}}

  \medskip 

  For distinction, the group we are going to construct for the proof of Theorem \ref{thm:main} will be denoted as $\Gamma '$ (instead of $\Gamma $ as in the case of the proof of Theorem \ref{thm:main2}). We will treat the case of $SL(2,\R)\times SL(2, \R)$ as a limit case of $SL(2,\R)\times SL(2, \C)$, more precisely, our 2-generated group $\Gamma '\leq SL(2,\R)\times SL(2,\R)$ will be a limit of 2-generated groups $\Gamma _n'\leq SL(2,\R)\times SL(2,\C)$.\footnote{This means that the generators of $\Gamma _n'$ converge to the corresponding generators of $\Gamma '$ in the $||\cdot ||$ norm.} 
  
  \medskip 
  
  Let us recall that in a real algebraic variety, the complement of the union of countably many subvarietes of positive co-dimension is dense. Then, for a dense subset $\mathcal{D}\subseteq SL(2,\R)\times SL(2, \R)$, any pair $(A,B)\in \mathcal{D}$ generates a non-Abelian free subgroup $\langle A, B\rangle $ of $SL(2,\R)$. Also, since $Q$ is hyperbolic, for any non-trivial word $w(X,Y)$, the relation $W(Q, X) = 1$ also defines a subvariety of $SL(2,\R)$ of a positive co-dimension. Then, for a dense subset $\mathcal{D}_0\subseteq SL(2,\R)$ and for any $A\in \mathcal{D}_0$, the pair $(Q,A)$ generates a non-Abelian free subgroup.  

  \medskip 

  We will use the matrix $Q$ from the previous section, but instead of $P$, we will work with a sequence of matrices $(P_n)$ in $SL(2,\C )$ satisfying certain properties as described below. 
  
\medskip 

  For all $n\geq 1$, let $$P_n =   \begin{bmatrix} 
	x_{11}^{(n)}   & x_{12}^{(n)}  \\
	x_{21}^{(n)}   & x_{22}^{(n)} 
 \end{bmatrix}$$

  where $$x_{ij}^{(n)} = p_{ij}^{(n)} + q_{ij}^{(n)}\beta  + r_{ij}^{(n)}\beta ^2  + s_{ij}^{(n)}\beta ^3$$ with $p_{ij}^{(n)}, q_{ij}^{(n)}, r_{ij}^{(n)}, s_{ij}^{(n)}\in \Z$ such that for all $i, j\in \{1,2\}$

   i)  $\displaystyle \mathop{\lim }_{n} (p_{ij}^{(n)} - r_{ij}^{(n)}\beta ^2) = u_{ij}$; 

   ii) $\displaystyle \mathop{\lim }_{n}(q_{ij}^{(n)}\beta  - s_{ij}^{(n)}\beta ^3) = 0$; 

   (iii) $\displaystyle \mathop{\lim }_{n} ((p_{11}^{(n)} + r_{11}^{(n)}\beta ^2) - (q_{11}^{(n)}\beta + s_{11}^{(n)}\beta ^3)) = v_{ij}$

   and the following conditions hold: 

\medskip 

   (iv) For all $n\geq 1$, the matrix $$R_n^{(1)} =  \begin{bmatrix} 
	(p_{11}^{(n)} + r_{11}^{(n)}\beta ^2) - (q_{11}^{(n)}\beta + s_{11}^{(n)}\beta ^3)   & (p_{12}^{(n)} + r_{12}^{(n)}\beta ^2) - (q_{12}^{(n)}\beta + s_{12}^{(n)}\beta ^3) \\
	(p_{21}^{(n)} + r_{21}^{(n)}\beta ^2) - (q_{21}^{(n)}\beta + s_{21}^{(n)}\beta ^3)  & (p_{22}^{(n)} + r_{22}^{(n)}\beta ^2) - (q_{22}^{(n)}\beta + s_{22}^{(n)}\beta ^3)
 \end{bmatrix}$$ is elliptic and the matrices $$R_n^{(2)} =  \begin{bmatrix} 
	(p_{11}^{(n)} - r_{11}^{(n)}\beta ^2) - (q_{11}^{(n)}\beta - s_{11}^{(n)}\beta ^3){\bf i}   & (p_{12}^{(n)} - r_{12}^{(n)}\beta ^2) - (q_{12}^{(n)}\beta - s_{12}^{(n)}\beta ^3){\bf i} \\
	(p_{21}^{(n)} - r_{21}^{(n)}\beta ^2) - (q_{21}^{(n)}\beta - s_{21}^{(n)}\beta ^3){\bf i}  & (p_{22}^{(n)} - r_{22}^{(n)}\beta ^2) - (q_{22}^{(n)}\beta - s_{22}^{(n)}\beta ^3){\bf i}
 \end{bmatrix}$$ and $$R_n^{(3)} =  \begin{bmatrix} 
	(p_{11}^{(n)} + r_{11}^{(n)}\beta ^2) + (q_{11}^{(n)}\beta + s_{11}^{(n)}\beta ^3)   & (p_{12}^{(n)} + r_{12}^{(n)}\beta ^2) + (q_{12}^{(n)}\beta + s_{12}^{(n)}\beta ^3) \\
	(p_{21}^{(n)} + r_{21}^{(n)}\beta ^2) + (q_{21}^{(n)}\beta + s_{21}^{(n)}\beta ^3)  & (p_{22}^{(n)} + r_{22}^{(n)}\beta ^2) + (q_{22}^{(n)}\beta + s_{22}^{(n)}\beta ^3)
 \end{bmatrix}$$ are hyperbolic;   

\medskip 

   (v) The matrix $R^{(1)} =  \begin{bmatrix} 
	v_{11}   & v_{12} \\
	v_{21}  & v_{22} 
 \end{bmatrix}$ is elliptic and the matrix $R^{(2)} =  \begin{bmatrix} 
	u_{11}   & u_{12} \\
	u_{21}  & u_{22} 
 \end{bmatrix}$ is hyperbolic; 

\medskip 

  (vi) The matrices $R^{(1)}$ and $Q$ generate a non-Abelian free group;

\medskip 

  (vii) $[QR^{(1)}Q^{-1}, R^{(1)}] \neq 1$ and $[\sigma _1(Q)R^{(2)}\sigma _1(Q)^{-1}, R^{(2)}] \neq 1$; 

  \medskip 

  (viii) for all $n\geq 1$, the matrices $R_n^{(3)}$ and $Q$ do not have a common eigenvector.

\medskip 

 Let $\Gamma _n' = \langle (Q, \sigma _1(Q), (R_n^{(1)}, R_n^{(2)}) \rangle , n\geq 1$ and for all natural $N\geq 1$, let $$\Gamma '(N) := \langle (Q^N, \sigma _1(Q)^N), ((R^{(1)})^N, (R^{(2)})^N) \rangle.$$ We have $\displaystyle \mathop{\lim }_nR_n^{(1)} = R^{(1)}$ and $\displaystyle \mathop{\lim }_nR_n^{(2)} = R^{(2)}$. We also note that the limit $\displaystyle \mathop{\lim }_nR_n^{(3)}$ does not necessarily exist as the entries may escape to infinity. 

 \medskip

  From the above conditions, as in the proof of Theorem \ref{thm:main2}, we obtain that the projections of $\Gamma '(N)$ onto both factors generate a dense subgroup in those factors. Thus, it remains to show the discreteness. For this, in addition to conditions (i)-(vii), we can also assume that for all $D > 0$, if $\epsilon >0$ is sufficiently small, there exists a natural $N$ such that for all $m\geq N, j\geq 1$ and for any non-identity word $W$, if $W(Q^{jm}, (R_n^{(3)})^{jm}) = \begin{bmatrix} 
	w_{11}   & w_{12} \\
	w_{21}  & w_{22} 
 \end{bmatrix}$ with $||W|| < \epsilon $, then for an eigenvector $[1, \lambda ]$ of $W$, we have  $$\mathrm{dist}([w_{11}:w_{21}], [1:\lambda ]) < D \ \mathrm{and} \ \mathrm{dist}([w_{21}:w_{22}], [1:\lambda ]) < D \ (3) $$.

  \medskip 
  
  Let us notice that in the proof of Theorem \ref{thm:main2}, in verifying the discreteness of $\Gamma $ (in the previous section), for sufficiently small $\epsilon > 0$ and $D > 0$, the choice of $N$ depends on $\epsilon $ and $D$ because we need to satisfy inequality (3) and also generate a non-Abelian free group in the third factor $\sigma _2(\Gamma _N)$. As $n\to \infty $, we can choose uniform $\epsilon $ and $D$ (i.e. both of these positive constants staying away from zero) to satisfy (3). On the other hand, as $n\to \infty $, the entries of the matrix $R_n^{(3)}$ can change erratically, however, if we already have a non-Abelian free group in the first factor (condition (vi)), the choice of $N$ would again be uniform for a sufficiently small $\epsilon $ and $D$. Thus we can claim that for some sufficiently small $\epsilon $ and $D$, and for sufficiently large $N$, for all sufficiently large $n$, the groups $\Gamma _{n, N}' = \langle (Q^N, \sigma _1(Q)^N), ((R_n^{(1)})^N, (R_n^{(2)})^N) \rangle$ are discrete in $SL(2,\R)\times SL(2,\C)$ with no non-identity element in the $\epsilon $-neighborhood of identity. Then the limit group $$\Gamma ' := \langle (Q^N, \sigma _1(Q)^N), ((R^{(1)})^N, (R^{(2)})^N) \rangle$$ is also discrete in $SL(2,\R)\times SL(2,\R)$.

 \end{document}